\documentclass[12pt]{article}

\newcommand{\qed}{{\hfill\rule{4pt}{7pt}}\medskip}

\usepackage{amsmath, epsfig}
\usepackage{amssymb}
\usepackage{amsfonts}
\usepackage{latexsym}

\newtheorem{thm}{Theorem}[section]
\newtheorem{prop}[thm]{Proposition}

\newtheorem{cor}[thm]{Corollary}

\newtheorem{lem}[thm]{Lemma}
\newtheorem{conj}[thm]{Conjecture}

\def\inv{{\rm inv}}

\def\pf{\noindent {\it Proof.} }

\numberwithin{equation}{section}

\begin{document}
\begin{center}
{\Large\bf Some Arithmetic Properties of the $q$-Euler Numbers\\[5pt]
 and  $q$-Sali\'e Numbers\footnote{European J. Combin. 27 (2006), 884--895.} }
\end{center}
\vskip 2mm
\centerline{Victor J. W. Guo$^1$  and Jiang Zeng$^2$}

\begin{center}
Institut Camille Jordan,
Universit\'e Claude Bernard (Lyon I)\\
F-69622, Villeurbanne Cedex, France \\
{\tt $^1$jwguo@eyou.com, $^2$zeng@igd.univ-lyon1.fr}
\end{center}

\vskip 0.7cm \noindent{\bf Abstract.}
For $m>n\geq 0$ and $1\leq d\leq m$, it is shown that the $q$-Euler number 
$E_{2m}(q)$ is congruent to $q^{m-n}E_{2n}(q)\bmod (1+q^d)$ if and only if 
$m\equiv n\bmod d$. The $q$-Sali\'e number $S_{2n}(q)$ is shown to be divisible by
$(1+q^{2r+1})^{\left\lfloor \frac{n}{2r+1}\right\rfloor}$ for any
$r\geq 0$. Furthermore, similar congruences for the generalized $q$-Euler numbers
are also obtained, and some conjectures are formulated.




\section{Introduction}
The Euler numbers $E_{2n}$ may be defined as the coefficients
in the Taylor expansion of $2/(e^x+e^{-x})$:
$$
\sum_{n=0}^{\infty}E_{2n}\frac{x^{2n}}{(2n)!}
=\left(\sum_{n=0}^\infty\frac{x^{2n}}{(2n)!}\right)^{-1}.
$$
A classical result due to Stern~\cite{Stern} asserts that
$$
E_{2m}\equiv E_{2n} \pmod{2^s} \quad \text{if and only if}
\quad 2m\equiv 2n \pmod{2^s}.
$$
The so-called Sali\'e numbers $S_{2n}$~\cite[p. 242]{GJ} are defined as
\begin{equation}
\sum_{n=0}^{\infty}S_{2n}\frac{x^{2n}}{(2n)!}=\frac{\cosh x}{\cos x}.
\label{eq:salie}
\end{equation}
Carlitz \cite{Carlitz65} first proved that 
the Sali\'e numbers $S_{2n}$ are divisible by $2^n$.

Motivated by the work of Andrews-Gessel \cite{AG}, Andrews-Foata \cite{AF},
D\'esarm\'enien \cite{De}, and Foata \cite{Foata}, we are about to study
a $q$-analogue of Stern's result and
a $q$-analogue of Carlitz's result for Sali\'e numbers.
A natural $q$-analogue of the Euler numbers is given by
\begin{align}
\sum_{n=0}^{\infty}E_{2n}(q)\frac{x^{2n}}{(q;q)_{2n}}
=\left(\sum_{n=0}^{\infty}\frac{x^{2n}}{(q;q)_{2n}}\right)^{-1},
\label{eq:qsecant}
\end{align}
where $(a;q)_{n}=(1-a)(1-aq)\cdots (1-aq^{n-1})$ for $n\geq 1$ and $(a;q)_0=1$.

A recent arithmetic study of
Euler numbers and more general $q$-Euler numbers can be found in \cite{Sun} and \cite{SZ}.
Note that, in order to coincide with
the Euler numbers in \cite{Sun,Wagstaff}, our definition of $E_{2n}(q)$
differs by a factor $(-1)^n$ from that in \cite{AF,AG,De,Foata}.

\begin{thm}\label{thm:e2m2n}
Let $m>n\geq 0$ and $1\leq d\leq m$. Then
$$
E_{2m}(q)\equiv q^{m-n}E_{2n}(q)
\pmod{1+q^d}
\quad\text{if and only if}\quad
m\equiv n \pmod d.
$$
\end{thm}

Since the polynomials $1+q^{2^ad}$ and $1+q^{2^b d}$ ($a\neq b$) are relatively prime, we
derive immediately from the above theorem the following
\begin{cor}
Let $m>n\geq 0$ and $2m-2n=2^s r$  with $r$ odd. Then
\begin{align*}
E_{2m}(q)\equiv q^{m-n}E_{2n}(q)
\pmod{\prod_{k=0}^{s-1}(1+q^{2^{k}r})}.
\end{align*}
\end{cor}

Define the $q$-Sali\'e numbers by
\begin{align}
\sum_{n=0}^{\infty}S_{2n}(q)\frac{x^{2n}}{(q;q)_{2n}}
=\sum_{n=0}^{\infty}\frac{q^nx^{2n}}{(q;q)_{2n}}
\left/\sum_{n=0}^{\infty}(-1)^n \frac{x^{2n}}{(q;q)_{2n}}.\right.
\label{eq:qsalie}
\end{align}
For each positive integer $n$, write $n=2^s(2r+1)$  with $r,s\geq 0$ (so $s$ is the
2-adic valuation of $n$), and  set $p_n(q)=1+q^{2r+1}$. Define
\begin{align*}
P_n(q)=\prod_{k=1}^{n}p_k(q)=\prod_{r\geq 0}(1+q^{2r+1})^{a_{n,r}},
\end{align*}
where $a_{n,r}$ is the number of
positive integers of the form $2^s(2r+1)$ less than or equal to $n$. The first values of
$P_{n}(q)$ are given in Table \ref{table:anr}.
\begin{table}[h]
{\footnotesize
\caption{Table of $P_{n}(q)$.\label{table:anr}}
\begin{center}
\begin{tabular}{|l|c|c|c|c|}
\hline
$n$             & 1 &  3 & 5  & 7   \\\hline
$P_{n}(q)$      & $(1+q)$ & $(1+q)^2(1+q^3)$  & $(1+q)^3(1+q^3)(1+q^5)$
  & $(1+q)^3(1+q^3)^2(1+q^5)(1+q^7)$ \\\hline
$n$             & 2 &  4 & 6  & 8   \\\hline
$P_{n}(q)$      & $(1+q)^2$  & $(1+q)^3(1+q^3)$ & $(1+q)^3(1+q^3)^2(1+q^5)$
  &  $(1+q)^4(1+q^3)^2(1+q^5)(1+q^7)$ \\\hline
\end{tabular}
\end{center}
}
\end{table}

Note that $P_{n}(1)=2^n$. The following is a $q$-analogue of Carlitz's result for
Sali\'e numbers:
\begin{thm}\label{thm:div-p}
For every $n\geq 1$, the polynomial $S_{2n}(q)$ is divisible by $P_n(q)$.
In particular, $S_{2n}(q)$ is divisible by
$(1+q^{2r+1})^{\left\lfloor \frac{n}{2r+1}\right\rfloor}$ for any
$r\geq 0$.
\end{thm}

We shall collect  some arithmetic properties of Gaussian
 polynomials or $q$-binomial coefficients in the next section.
The proofs of Theorems \ref{thm:e2m2n} and \ref{thm:div-p} are given in Sections 3 and 4, respectively.
We will give some similar arithmetic properties of the generalized $q$-Euler
numbers in Section 5.
Some combinatorial remarks and open problems are given in Section 6.

\section{Two properties of Gaussian polynomials}
The Gaussian polynomial ${M\brack N}_q$ may be defined by
\[
{M\brack N}_q=
\begin{cases}
\displaystyle\frac{(q;q)_M}{(q;q)_N (q;q)_{M-N}}, &\text{if $0\leq N\leq M$},\\[5pt]
0,&\text{otherwise.}
\end{cases}
\]
The following result is equivalent to the so-called $q$-Lucas theorem (see 
Olive \cite{Olive} and D\'esarm\'enien \cite[Proposition 2.2]{De}).
\begin{prop}\label{prop:root}
Let $m,k,d$ be positive integers, and write $m=ad+b$ and $k=rd+s$,
where $0\leq b,s\leq d-1$. Let $\omega$ be a primitive $d$-th root
of unity.  Then
$$
{m\brack k}_{\omega}={a\choose r}{b\brack s}_\omega.
$$
\end{prop}
Indeed, we have
\begin{align*}
{m\brack k}_q
&=\prod_{j=1}^{rd+s}\frac{1-q^{(a-r)d+b-s+j}}{1-q^j}\\[5pt]
&=\left(\prod_{j=1}^{s}\frac{1-q^{(a-r)d+b-s+j}}{1-q^j}\right)
\left(\prod_{j=1}^{rd}\frac{1-q^{(a-r)d+b+j}}{1-q^{s+j}}\right).
\end{align*}
By definition, we have $\omega^d=1$ and $\omega^j\neq 1$
for $0<j<d$. Hence,
\begin{align*}
\lim_{q\rightarrow \omega}\prod_{j=1}^{s}\frac{1-q^{(a-r)d+b-s+j}}{1-q^j}
&=\prod_{j=1}^{s}\frac{1-\omega^{b-s+j}}{1-\omega^j}={b\brack s}_\omega.
\end{align*}
Notice that, for any integer $k$, the set $\{k+j\colon j=1,\ldots, rd\}$ is
a {\it complete system of residues} modulo $rd$. Therefore,
\begin{align*}
\lim_{q\rightarrow \omega}\prod_{j=1}^{rd}\frac{1-q^{(a-r)d+b+j}}{1-q^{s+j}}
&=\lim_{q\rightarrow \omega}\frac{(1-q^{(a-r+1)d})(1-q^{(a-r+2)d})\cdots (1-q^{ad})}
{(1-q^d)(1-q^{2d})\cdots(1-q^{rd})}\\[5pt]
&={a\choose r}.
\end{align*}

Let $\Phi_n(x)$ be the $n$-th {\it cyclotomic polynomial}.
The following easily proved result can be found in \cite[Equation (10)]{KW}.
\begin{prop}\label{prop:factor}
The Gaussian polynomial ${m\brack k}_q$ can be factorized into
$$
{m\brack k}_q=\prod_{d}\Phi_d(q),
$$
where the product is over all positive integers $d\leq m$ such that
$\lfloor k/d\rfloor+\lfloor (m-k)/d\rfloor<\lfloor m/d\rfloor$.
\end{prop}
Indeed, using the factorization
$q^n-1=\prod_{d|n}\Phi_d(q)$,
we have
$$
(q;q)_m=(-1)^m \prod_{k=1}^{m}\prod_{d\mid k}\Phi_d(q)
=(-1)^m \prod_{d=1}^{m}\Phi_d(q)^{\lfloor m/d\rfloor},
$$
and so
$$
{m\brack k}_q=\frac{(q;q)_m}{(q;q)_k (q;q)_{m-k}}
=\prod_{d=1}^{m}\Phi_d(q)^{\lfloor m/d\rfloor
-\lfloor k/d\rfloor-\lfloor (m-k)/d\rfloor}.
$$
Proposition \ref{prop:factor} now follows from the obvious fact that
\begin{align*}
\lfloor \alpha+\beta\rfloor-\lfloor \alpha\rfloor -\lfloor \beta\rfloor
=0\quad \textrm{or}\quad 1,\quad\text{for $\alpha,\beta\in\mathbb{R}$.}
\end{align*}

\section{Proof of Theorem \ref{thm:e2m2n}}
Multiplying both sides of \eqref{eq:qsecant} by
$\sum_{n=0}^{\infty}x^{2n}/(q;q)_{2n}$ and equating coefficients of $x^{2m}$,
we see that $E_{2m}(q)$ satisfies the following recurrence relation:
\begin{align}
E_{2m}(q)=-\sum_{k=0}^{m-1}{2m\brack 2k}_q E_{2k}(q).\label{eq:rece2n}
\end{align}
This enables us to obtain the first values of the $q$-Euler numbers:
\begin{align*}
E_0(q)&=-E_2(q)=1,\\[5pt]
E_4(q)&=q(1+q)(1+q^2)+q^2,\\[5pt]
E_6(q)&=-q^2(1+q^3)(1+4q+5q^2+7q^3+6q^4+5q^5+2q^6+q^7)+q^3.
\end{align*}
We first establish the following result.
\begin{lem}\label{lem:relem}
Let $m>n\geq 0$ and $1\leq d\leq m$. Then
\begin{equation}\label{eq:relem}
E_{2m}(q)\equiv q^{m-n}E_{2n}(q) \pmod{\Phi_{2d}(q)}
\text{\quad if and only if\quad}
m\equiv n \pmod d.
\end{equation}
\end{lem}
\pf It is easy to see that Lemma \ref{lem:relem} is equivalent to
\begin{equation}\label{eq:zeta}
E_{2m}(\zeta)=\zeta^{m-n}E_{2n}(\zeta)
\quad \text{if and only if}
\quad m\equiv n \pmod d,
\end{equation}
where $\zeta\in\mathbb{C}$ is a $2d$-th primitive root of unity. 

We proceed by induction on $m$. Statement \eqref{eq:zeta} is trivial for $m=1$.
Suppose it holds for every number less than $m$.
Let $n<m$ be fixed. Write $m=ad+b$ with $0\leq b\leq d-1$, then $2m=a(2d)+2b$. 
By Proposition~\ref{prop:root}, we see that
\begin{equation}\label{eq:2m2kzeta}
{2m\brack 2k}_{\zeta}
={a\choose r}{2b\brack 2s}_\zeta,\quad\text{where}\quad k=rd+s,\ 0\leq s\leq d-1.
\end{equation}
Hence, by \eqref{eq:rece2n} and \eqref{eq:2m2kzeta}, we have
\begin{align}
E_{2m}(\zeta)
&=-\sum_{k=0}^{m-1}{2m\brack 2k}_{\zeta}E_{2k}(\zeta) \nonumber\\[5pt]
&=-\sum_{r=0}^{a}\sum_{s=0}^{b-\delta_{a\,r}}
{a\choose r}{2b\brack 2s}_\zeta E_{2rd+2s}(\zeta), \nonumber\\[5pt]
&=-\sum_{s=0}^{b}\sum_{r=0}^{a-\delta_{b\,s}}
{a\choose r}{2b\brack 2s}_\zeta E_{2rd+2s}(\zeta), \label{eq:2mto}
\end{align}
where $\delta_{i\,j}$ equals $1$ if $i=j$ and $0$ otherwise.

By the induction hypothesis, we have
\begin{align}
E_{2rd+2s}(\zeta)=\zeta^{rd}E_{2s}(\zeta)=(-1)^{r}E_{2s}(\zeta).
\label{eq:kizeta}
\end{align}
Thus, 
\begin{align*}
\sum_{r=0}^{a}
{a\choose r}{2b\brack 2s}_\zeta E_{2rd+2s}(\zeta)
={2b\brack 2s}_\zeta E_{2s}(\zeta)\sum_{r=0}^{a}
{a\choose r}(-1)^{r}=0.
\end{align*}
Therefore, Equation \eqref{eq:2mto} implies that
\begin{equation}\label{eq:e2mz}
E_{2m}(\zeta)=(-1)^{a}E_{2b}(\zeta)=\zeta^{m-b}E_{2b}(\zeta).
\end{equation}
{}From \eqref{eq:e2mz} we see that
$$
E_{2m}(\zeta)=\zeta^{m-n}E_{2n}(\zeta)\quad
\Longleftrightarrow \quad
E_{2n}(\zeta)=\zeta^{n-b}E_{2b}(\zeta).
$$
By the induction hypothesis, the latter equality is also equivalent to
$$n\equiv b \pmod d \quad\Longleftrightarrow\quad m\equiv n \pmod d.$$
This completes the proof. \qed

Since
$$
1+q^d={\frac{q^{2d}-1}{q^{d}-1}}
=\frac{\prod_{k\mid 2d}\Phi_k(q)}{\prod_{k\mid d}\Phi_k(q)}
=\prod_{\substack{k\mid d \\2k\nmid d}}\Phi_{2k}(q),
$$
and any two different cyclotomic polynomials are relatively prime, 
Theorem \ref{thm:e2m2n} follows from Lemma~\ref{lem:relem}.

\medskip
\noindent{\it Remark.} The sufficiency part of \eqref{eq:relem} is equivalent to
D\'esarm\'enien's result \cite{De}:
$$
E_{2km+2n}(q)\equiv (-1)^{m}E_{2n}(q)
\pmod{\Phi_{2k}(q)}.
$$

\section{Proof of Theorem \ref{thm:div-p}}
Recall that the $q$-tangent numbers $T_{2n+1}(q)$ are defined by
\[
\sum_{n=0}^{\infty}T_{2n+1}(q)\frac{x^{2n+1}}{(q;q)_{2n+1}}
=\sum_{n=0}^{\infty}(-1)^n\frac{x^{2n+1}}{(q;q)_{2n+1}}
\left/\sum_{n=0}^{\infty}(-1)^n
\frac{x^{2n}}{(q;q)_{2n}}\right..
\]
Foata~\cite{Foata} proved that $T_{2n+1}(q)$ is divisible by $D_n(q)$, where
\begin{align*}
D_n(q)=
\begin{cases}
\displaystyle\prod_{k=1}^{n}Ev_k(q), &\text{if $n$ is odd,}\\[10pt]
\displaystyle(1+q^2)\prod_{k=1}^{n}Ev_k(q), &\text{if $n$ is even,}
\end{cases}
\end{align*}
and
\begin{align*}
Ev_n(q)=\prod_{j=0}^{s}(1+q^{2^j r}), \quad\text{where $n=2^s r$ with $r$ odd.}
\end{align*}

Notice that this implies
that $T_{2n+1}(q)$ is divisible by both
$(1+q)^n$ and $(-q;q)_n$, a result
due to Andrews and Gessel \cite{AG}.

To prove our theorem we need the following relation relating $S_{2n}(q)$ 
to $T_{2n+1}(q)$.
\begin{lem}\label{lem:s2n}
For every $n\geq 1$, we have
\begin{align}
\sum_{k=0}^{n}(-1)^k q^{k}{2n\brack 2k}_q S_{2k}(q)S_{2n-2k}(q)
=T_{2n-1}(q)(1-q^{2n}). \label{eq:inv-e2n}
\end{align}
\end{lem}
\pf 
Replacing $x$ by $q^{1/2}ix$ ($i=\sqrt{-1}$) in \eqref{eq:qsalie},
we obtain
\begin{align}
\sum_{n=0}^{\infty}S_{2n}(q)\frac{(-1)^nq^nx^{2n}}{(q;q)_{2n}}
=\sum_{n=0}^{\infty}\frac{(-1)^nq^{2n}x^{2n}}{(q;q)_{2n}}
\left/\sum_{n=0}^{\infty}\frac{q^nx^{2n}}{(q;q)_{2n}}.\right.
\label{eq:qisal}
\end{align}
Multiplying \eqref{eq:qsalie} with \eqref{eq:qisal}, we get
\begin{align}
&\left(\sum_{n=0}^{\infty}S_{2n}(q)\frac{x^{2n}}{(q;q)_{2n}}\right)
\left(\sum_{n=0}^{\infty}S_{2n}(q)\frac{(-1)^nq^nx^{2n}}{(q;q)_{2n}}\right)
\label{eq:qihalf}\\[5pt]
&=\sum_{n=0}^{\infty}(-1)^n\frac{q^{2n}x^{2n}}{(q;q)_{2n}}
\left/\sum_{n=0}^{\infty}(-1)^n\frac{x^{2n}}{(q;q)_{2n}}\right.
\nonumber\\[5pt]
&=1+x\sum_{n=1}^{\infty}(-1)^{n-1}\frac{q^{2n-1}x^{2n-1}}{(q;q)_{2n-1}}
\left/\sum_{n=0}^{\infty}(-1)^n\frac{x^{2n}}{(q;q)_{2n}}\right.\nonumber\\[5pt]
&=1+x\sum_{n=0}^{\infty}T_{2n+1}(q)\frac{x^{2n+1}}{(q;q)_{2n+1}}. \label{eq:ix-ix}
\end{align}
Equating the coefficients of $x^{2n}$ in \eqref{eq:qihalf} and \eqref{eq:ix-ix},
we are led to \eqref{eq:inv-e2n}. \qed

\medskip

It is easily seen that $P_n(q)$ is the \emph{least common multiple} of
the polynomials $(1+q^{2r+1})^{\left\lfloor \frac{n}{2r+1}\right\rfloor}$
($r\geq 0$).
For any $r\geq 0$, there holds
$$
1+q^{2r+1}=\frac{q^{4r+2}-1}{q^{2r+1}-1}
=\frac{\prod_{d\mid (4r+2)}\Phi_d(q)}{\prod_{d\mid (2r+1)}\Phi_d(q)}
=\prod_{d\mid (2r+1)}\Phi_{2d}(q).
$$
It follows that
$$
P_n(q)=\prod_{r\geq 0}\Phi_{4r+2}(q)^{\lfloor \frac{n}{2r+1}\rfloor}.
$$
Theorem \ref{thm:div-p} is trivial for $n=1$.
Suppose it holds for all integers less than $n$.
In the summation of the left-hand side of \eqref{eq:inv-e2n},
combining the first and last terms, we can rewrite
 Equation~\eqref{eq:inv-e2n} as follows:
\begin{equation}\label{eq:casesn}
(1+(-1)^n q^n)S_{2n}(q)
+\sum_{k=1}^{n-1}(-1)^k q^{k}{2n\brack 2k}_q S_{2k}(q)S_{2n-2k}(q)
=T_{2n-1}(q)(1-q^{2n}).
\end{equation}

For every $k$ ($1\leq k\leq n-1$),
by the induction hypothesis, the polynomial $S_{2k}(q)S_{2n-2k}(q)$ is divisible by
$$
P_{k}(q)P_{n-k}(q)
=\prod_{r\geq 0}\Phi_{4r+2}(q)^{\lfloor \frac{k}{2r+1}\rfloor
+\lfloor \frac{n-k}{2r+1}\rfloor}.
$$
And by Proposition~\ref{prop:factor}, we have
$$
{2n\brack 2k}_q
=\prod_{d=1}^{2n}\Phi_d(q)^{\lfloor 2n/d\rfloor
-\lfloor 2k/d\rfloor-\lfloor (2n-2k)/d\rfloor},
$$
which is clearly divisible by
$$
\prod_{r\geq 0}\Phi_{4r+2}(q)^{
\lfloor \frac{n}{2r+1}\rfloor
-\lfloor \frac{k}{2r+1}\rfloor
-\lfloor \frac{n-k}{2r+1}\rfloor}.
$$
Hence, the product ${2n\brack 2k}_q S_{2k}(q)S_{2n-2k}(q)$ is divisible by
$$
\prod_{r\geq 0}\Phi_{4r+2}(q)^{\lfloor \frac{n}{2r+1}\rfloor}=P_n(q).
$$
Note that $P_{n-1}(q)\mid D_{n-1}(q)$ and $p_n(q)\mid (1-q^{2n})$.
Therefore, by \eqref{eq:casesn} and the aforementioned result of Foata, 
we immediately have
$$
P_n(q)\mid (1+(-1)^nq^n)S_{2n}(q).
$$
Since $P_n(q)$ is relatively prime to $(1+(-1)^nq^n)$, we obtain $P_n(q)\mid S_{2n}(q)$. 

\medskip
\noindent{\it Remark.}
Since $S_0(q)=1$ and $S_2(q)=1+q$,
using \eqref{eq:casesn} and the divisibility of $T_{2n+1}(q)$, we can prove
by induction that $S_{2n}(q)$ is divisible by $(1+q)^n$ without
using the divisibility property of Gaussian polynomials.

\section{The generalized $q$-Euler numbers}
The generalized Euler numbers may be defined by
\begin{align*}
\sum_{n=0}^{\infty}E_{kn}^{(k)}\frac{x^{kn}}{(kn)!}
=\left(\sum_{n=0}^{\infty}\frac{x^{kn}}{(kn)!}\right)^{-1}.
\end{align*}
Some congruences for these numbers are given in \cite{Gessel,LM}.
A $q$-analogue of generalized Euler numbers is given by
\begin{align*}
\sum_{n=0}^{\infty}E_{kn}^{(k)}(q)\frac{x^{kn}}{(q;q)_{kn}}
=\left(\sum_{n=0}^{\infty}\frac{x^{kn}}{(q;q)_{kn}}\right)^{-1},
\end{align*}
or, recurrently,
\begin{align}
E_{0}^{(k)}(q)=1,\quad
E_{kn}^{(k)}(q)=-\sum_{j=0}^{n-1}{kn\brack kj}_q E_{kj}^{(k)}(q),
\quad n\geq 1.
\label{eq:gen-rec}
\end{align}
Note that $E_{kn}^{(k)}(q)$ is equal to $(-1)^n f_{nk,k}(q)$
studied by Stanley \cite[p.~148, Equation (57)]{Stanley97}. 
\begin{thm}\label{thm:genezeta}
Let $m>n\geq 0$ and $1\leq d\leq m$. Let $k\geq 1$, and let $\zeta\in\mathbb{C}$ 
be a $2kd$-th primitive root of unity. Then
\begin{align}
E_{km}^{(k)}(\zeta^2)=\zeta^{k(m-n)}E_{kn}^{(k)}(\zeta^2)\label{eq:genezeta}
\end{align}
if and only if
$$
m\equiv n \pmod d.
$$
\end{thm}
The proof is by induction on $m$ and using the recurrence
definition \eqref{eq:gen-rec}.
Since it is  analogous to the proof of \eqref{eq:zeta}, we omit it here.
Note that $\zeta^2$ in Theorem~\ref{thm:genezeta}
is a $kd$-th primitive root of unity. Therefore, when $k$ is even
or $m\equiv n \bmod 2$, Equation~\eqref{eq:genezeta} is equivalent to
$$
E_{km}^{(k)}(q)\equiv q^{\frac{k(m-n)}2}E_{kn}^{(k)}(q) \pmod{\Phi_{kd}(q)}.
$$

As mentioned before,
$$
1+q^{2^kd}=\prod_{\substack{i\mid 2^kd\\2i\,\nmid\, 2^kd}}\Phi_{2i}(q),
$$
and we obtain the following theorem and its corollaries.
\begin{thm}
Let $k\geq 1$. Let $m>n\geq 0$ and $1\leq d\leq m$. Then
$$
E_{2^km}^{(2^k)}(q)\equiv q^{2^{k-1}(m-n)}E_{2^kn}^{(2^k)}(q)
 \pmod{1+q^{2^{k-1}d}}
\quad\text{if and only if}\quad
m\equiv n \pmod d.
$$
\end{thm}

\begin{cor}
Let $k\geq 1$. Let $m>n\geq 0$ and $m-n=2^{s-1} r$ with $r$ odd. Then
\begin{align*}
E_{2^km}^{(2^k)}(q)\equiv q^{2^{k-1}(m-n)}E_{2^kn}^{(2^k)}(q)
\pmod{\prod_{i=0}^{s-1}(1+q^{2^{k+i-1}r})}.
\end{align*}
\end{cor}

\begin{cor}Let $k,m,n,s$ be as above. Then
\begin{align*}
E_{2^km}^{(2^k)}\equiv E_{2^kn}^{(2^k)} \pmod{2^s}.
\end{align*}
\end{cor}

Furthermore, numerical evidence seems to suggest the following congruence conjecture for
generalized Euler numbers.
\begin{conj}
Let $k\geq 1$. Let $m>n\geq 0$ and $m-n=2^{s-1} r$ with $r$ odd. Then
\begin{align*}
E_{2^km}^{(2^k)}\equiv E_{2^kn}^{(2^k)}+2^s \pmod{2^{s+1}}.
\end{align*}
\end{conj}
This conjecture is clearly a generalization of Stern's result, which corresponds to
the $k=1$ case.

\section{Concluding remarks}
We can also consider the following variants
of
the $q$-Sali\'e numbers:
\begin{align}
\sum_{n=0}^{\infty}\overline{S}_{2n}(q)\frac{x^{2n}}{(q;q)_{2n}}
&=\sum_{n=0}^{\infty}\frac{x^{2n}}{(q;q)_{2n}}\left/
\sum_{n=0}^{\infty}(-1)^n \frac{x^{2n}}{(q;q)_{2n}}\right.,
\label{eq:qsalie2}\\[5pt]
\sum_{n=0}^{\infty}\widehat{S}_{2n}(q)\frac{x^{2n}}{(q;q)_{2n}}
&=\sum_{n=0}^{\infty}\frac{q^{2n}x^{2n}}{(q;q)_{2n}}
\left/\sum_{n=0}^{\infty}(-1)^n \frac{x^{2n}}{(q;q)_{2n}}\right.,
\label{eq:qtsalie3}\\[5pt]
\sum_{n=0}^{\infty}\widetilde{S}_{2n}(q)\frac{x^{2n}}{(q;q)_{2n}}
&=\sum_{n=0}^{\infty}\frac{q^{n^2}x^{2n}}{(q;q)_{2n}}\left/
\sum_{n=0}^{\infty}(-1)^n \frac{x^{2n}}{(q;q)_{2n}}\right..
\label{eq:qtsalie}
\end{align}
Multiplying both sides of \eqref{eq:qsalie2}--\eqref{eq:qtsalie} by
$\sum_{n=0}^{\infty}(-1)^n x^{2n}/(q;q)_{2n}$ and equating
coefficients of $x^{2n}$, we obtain
\begin{align}
\overline{S}_{2n}(q)&=1-\sum_{k=0}^{n-1}(-1)^{n-k}{2n\brack 2k}_q
\overline{S}_{2k}(q), \label{eq:salone}\\[5pt]
\widehat{S}_{2n}(q)&=q^{2n}-\sum_{k=0}^{n-1}(-1)^{n-k}{2n\brack 2k}_q
\widehat{S}_{2k}(q),\label{eq:saltwo}\\[5pt]
\widetilde{S}_{2n}(q)&=q^{n^2}-\sum_{k=0}^{n-1}(-1)^{n-k}{2n\brack 2k}_q
\widetilde{S}_{2k}(q). \label{eq:salthree}
\end{align}
This gives
\begin{align*}
&\overline{S}_0(q)=1, \quad \overline{S}_2(q)=2,
\quad \overline{S}_4(q)=2(1+q^2)(1+q+q^2), \\[5pt]
&\widehat{S}_0(q)=1, \quad \widehat{S}_2(q)=1+q^2,
\quad \widehat{S}_4(q)=q(1+q^2)(1+3q+q^2+q^3), \\[5pt]
&\widetilde{S}_0(q)=1, \quad \widetilde{S}_2(q)=1+q,
\quad \widetilde{S}_4(q)=q(1+q)(1+q^2)(2+q), 
\end{align*}
and
\begin{align*}
&\overline{S}_6(q)=
2(1+q^2)(1+q+2q^2+4q^3+6q^4+6q^5+6q^6+5q^7+4q^8+2q^9+q^{10}), \\[5pt]
&\widehat{S}_6(q)=
q^{2}(1+q^2)^2(1+4q+7q^2+6q^3+6q^4+6q^5+5q^6+2q^7+q^8),\\[5pt]
&\widetilde{S}_6(q)=q^2(1+q)(1+q^2)(1+q^3)(2+4q+5q^2+4q^3+3q^4+q^5).
\end{align*}

For $n\geq 1$ define three sequences of polynomials:
\begin{align*}
\overline{Q}_n(q)&:=\prod_{r\geq 1}\Phi_{4r}(q)^{\lfloor
\frac{n}{2r}\rfloor},\\[5pt]
\widehat{Q}_n(q)&:=
\begin{cases}
\overline{Q}_n(q), &\text{if $n$ is even,}\\[5pt]
(1+q^2)\overline{Q}_n(q), &\text{if $n$ is odd,}
\end{cases}\\[5pt]
\widetilde{Q}_n(q)&:=(1+q)(1+q^2)\cdots (1+q^n).
\end{align*}
Note that $\overline{Q}_n(q)$ is the least common multiple of the polynomials
$(1+q^{2r})^{\lfloor \frac{n}{2r}\rfloor}$, $r\geq 1$ (see Table~\ref{table:qover}). 
\begin{table}[h]
{\footnotesize
\caption{Table of $\overline{Q}_n(q)$.\label{table:qover}}
\begin{center}
\begin{tabular}{|l|c|c|c|c|}
\hline
$n$             & 1 &  3 & 5  & 7   \\\hline
$\overline{Q}_n(q)$      & $1$ & $1+q^2$  & $(1+q^2)^2(1+q^4)$
& $(1+q^2)^2(1+q^4)(1+q^6)$ \\\hline
$n$             & 2 &  4 & 6  & 8   \\\hline
$\overline{Q}_n(q)$      & $1+q^2$  & $(1+q^2)^2(1+q^4)$ & $(1+q^2)^2(1+q^4)^2(1+q^6)$
  &  $(1+q^2)^3(1+q^4)^2(1+q^6)(1+q^8)$ \\\hline
\end{tabular}
\end{center}
}
\end{table}

From \eqref{eq:salone}--\eqref{eq:salthree}, it is easy to derive 
by induction that for $n\geq 1$,
$$
2 \mid  \overline{S}_{2n}(q),\quad
(1+q^2)  \mid \widehat{S}_{2n}(q),\quad
(1+q)\mid \widetilde{S}_{2n}(q).
$$
Moreover, the computation of the first values of these polynomials seems to
suggest the following stronger result.
\begin{conj}\label{conj:sal2}
For $n\geq 1$, we have the following divisibility properties:
$$
\overline{Q}_n(q) \mid  \overline{S}_{2n}(q),\quad
\widehat{Q}_n(q)  \mid \widehat{S}_{2n}(q),\quad
\widetilde{Q}_n(q)\mid \widetilde{S}_{2n}(q).
$$
\end{conj}

Similarly to the proof of Lemma~\ref{lem:s2n}, we can obtain
$$
\left(\sum_{n=0}^{\infty}\widehat{S}_{2n}(q)\frac{x^{2n}}{(q;q)_{2n}}\right)
\left(\sum_{n=0}^{\infty}\widehat{S}_{2n}(q)\frac{(-1)^nq^{2n}x^{2n}}{(q;q)_{2n}}\right)
=1-qx^2+(1+q)x\sum_{n=0}^{\infty}T_{2n+1}(q)\frac{x^{2n+1}}{(q;q)_{2n+1}},
$$
which yields
\begin{align}
\sum_{k=0}^{n}(-1)^k q^{2k}{2n\brack 2k}_q \widehat{S}_{2k}(q)\widehat{S}_{2n-2k}(q)
=T_{2n-1}(q)(1+q)(1-q^{2n}), \quad n\geq 2. \label{eq:inv-sal3}
\end{align}
However, it seems difficult to use \eqref{eq:inv-sal3} to prove directly
the divisibility of $\widehat{S}_{2n}(q)$ by $\widehat{Q}_n(q)$,
because when $n$ is even $1+(-1)^nq^{2n}$ is in general not
relatively prime to $\widehat{Q}_n(q)$.

Finally it is well-known that $E_{2n}(q)$ has a nice combinatorial
interpretation in terms of generating functions of alternating
permutations. Recall that a \emph{permutation} $x_1 x_2 \cdots
x_{2n}$ of $[2n]:=\{1,2,\ldots, 2n\}$ is called {\it alternating},
if $x_1<x_2>x_3<\cdots>x_{2n-1}<x_{2n}$. As usual, the number of
\emph{inversions} of a permutation $x=x_1 x_2\cdots x_n$, denoted
$\inv(x)$, is defined to the number of pairs $(i,j)$ such that
$i<j$ and $x_i>x_j$. It is known (see \cite[p.~148, Proposition
3.16.4]{Stanley97}) that
$$
(-1)^nE_{2n}(q)=\sum_{\pi}q^{\inv(\pi)},
$$
where $\pi$ ranges over all the alternating permutations of
$[2n]$. It would be interesting to find a combinatorial proof of
Theorem 1 within the alternating permutations model.

A permutation $x=x_1 x_2 \cdots x_{2n}$ of $[2n]$ is said to be
a {\it Sali\'e permutation}, if there exists an even index $2k$
such that $x_1x_2\cdots x_{2k}$ is alternating and
$x_{2k}<x_{2k+1}<\cdots <x_{2n}$, and $x_{2k-1}$ is called the
{\it last valley} of $x$. It is known (see \cite[p.~242, Exercise
4.2.13]{GJ}) that $\frac{1}{2}S_{2n}$ is the number of Sali\'e permutations of $[2n]$.

\begin{prop}
For every $n\geq 1$ the polynomial
$\frac{1}{2}\overline{S}_{2n}(q)$ is the generating function for
Sali\'e permutations of $[2n]$ by number of inversions.
\end{prop}
\pf
Substituting \eqref{eq:qsecant} into \eqref{eq:qsalie2}
and comparing coefficients of $x^{2n}$ on both sides, we obtain
\begin{equation}\label{eq:sal-eul}
\overline{S}_{2n}(q)=\sum_{k=0}^{n} {2n\brack 2k}_q (-1)^k E_{2k}(q).
\end{equation}
As ${2n\brack 2k}_q$ is the generating function for the
permutations of $1^{2k}2^{2n-2k}$ by number of inversions (see
e.g. \cite[p.~26, Proposition~1.3.17]{Stanley97}), it is easily seen that
${2n\brack 2k}_q (-1)^kE_{2k}(q)$ is the generating function for
permutations $x=x_1x_2\cdots x_{2n}$ of $[2n]$ such that $x_1x_2\cdots x_{2k}$
is alternating and $x_{2k+1}\cdots x_{2n}$ is increasing with respect to number of inversions.
Notice that such a permutation $x$
is a Sali\'e permutation with the last valley $x_{2k-1}$ if $x_{2k}<x_{2k+1}$ or
$x_{2k+1}$ if $x_{2k}>x_{2k+1}$. Therefore, the right-hand side of
\eqref{eq:sal-eul} is twice the generating function for Sali\'e
permutations of $[2n]$ by number of inversions. This completes the
proof. \qed

It is also possible to find similar combinatorial interpretations for
the other $q$-Sali\'e numbers, which are left to the interested readers.


\vskip 5mm
\noindent{\bf Acknowledgment.}
The second author was supported by EC's IHRP Programme, within Research Training
Network ``Algebraic Combinatorics in Europe," grant HPRN-CT-2001-00272.

\renewcommand{\baselinestretch}{1}

\end{document}